# APPROXIMATE SOLUTIONS OF COUPLED RAMANI EQUATION BY USING RDTM WITH COMPARED DTM AND EXACT SOLUTIONS

Murat GUBES[1], Galip OTURANC[2]

*In this paper, we present the new approximate solutions of famous coupled Ramani Equation. In order to obtain the solution, we use the semi-analytic methods differential transform method (DTM) and reduced form of DTM called reduced differential transform method (RDTM). We compare our proposed methods with exact solution and also compare with together. Numerical results show clearly that DTM and RDTM are very effective and also provide very accurate solutions. Additionally, it should be noted that RDTM is applied very easily, fast and more convergent than DTM for these kind of problems.*

## 1. Introduction

Partial differential equations are the fundamental phenomena and applied area in physic, engineering, chemistry and etc. In real life, many events can be modeled by a nonlinear partial differential equation such as evolution equations. Particularly in nonlinear sciences, one of the important and outstanding evolution equation is the famous coupled Ramani Equation that is presented as follow [1,2,3,4]

$$u(x,t)_{xxxxxx} + 15u(x,t)_{xx} u(x,t)_{xxx} + 15u(x,t)_x u(x,t)_{xxxx}$$
$$+ 45u(x,t)_x^2 u(x,t)_{xx} - 5u(x,t)_{tt} + 18v(x,t)_x$$
$$- 5\left[ u(x,t)_{xxxt} + 3u(x,t)_{xx} u(x,t)_t + 3u(x,t)_x u(x,t)_{xt} \right] = 0 \qquad (1)$$

$$v(x,t)_t - v(x,t)_{xxx} - 3v(x,t)_x u(x,t)_x - 3v(x,t)u(x,t)_{xx} = 0$$

---

[1] Dept.of Mathematics, Karamanoglu Mehmetbey University, Karaman, TURKEY,
e-mail: mgubes@kmu.edu.tr
[2] Dept.of Mathematics, Selcuk University, Konya, TURKEY, e-mail:goturanc@sekcuk.edu.tr

In literature, a great number of researchers have studied the system (1) to obtain exact and approximate solutions as seen some of in [5-16]. Ablowitz and Clarkson [5], Ito [6], Zhang [7], Feng [8], Malfiet and Hereman [9] have investigated the solitons and inverse scattering, extensions, exact traveling wave solutions and traveling solitary wave solutions of nonlinear evolution equations respectively. Li has presented exact traveling wave solutions of six order Ramani and a coupled Ramani equation in [10]. In [11], Nadjafikhah and Shirvani-Sh have found Lie symmetries and conservation laws of Hirota-Ramani equation. Further, Yusufoğlu and Bekir have obtained the two exact traveling wave solutions of coupled Ramani equation by applying *tanh method* as following [12]

$$u(x,t) = a_0 + 2\alpha \tanh[\alpha(x-\beta t)],$$
$$v(x,t) = -\frac{4}{9}\beta\alpha^4 - \frac{16}{27}\alpha^6 + \frac{5}{9}\alpha^2\beta^2 - \frac{5}{54}\beta^3$$
$$+ \left(\frac{20}{9}\beta\alpha^4 + \frac{16}{9}\alpha^6 - \frac{5}{9}\alpha^2\beta^2\right)\tanh^2[\alpha(x-\beta t)] \quad (2)$$

and

$$u(x,t) = a_0 - 2\alpha \tan[\alpha(x-\beta t)],$$
$$v(x,t) = -\frac{4}{9}\beta\alpha^4 + \frac{16}{27}\alpha^6 - \frac{5}{9}\alpha^2\beta^2 - \frac{5}{54}\beta^3 \quad, \left(|\alpha(x-\beta t)| < \frac{\pi}{2}\right) \quad (3)$$
$$+ \left(-\frac{20}{9}\beta\alpha^4 + \frac{16}{9}\alpha^6 - \frac{5}{9}\alpha^2\beta^2\right)\tan^2[\alpha(x-\beta t)]$$

where $a_0, \alpha$ and $\beta$ are arbitrary constants.

Recently, Wazwaz and Triki [13], Wazwaz [14], Jafarian et all [15] and Wazwaz [16] have presented the multiple soliton solutions and approximate solution of eq. (1) respectively.

According to the this, main aim of our study is to obtain accurate, convergent and efficient approximate solution of coupled Ramani equation (1) by using differential transform (DTM) and reduced differential transform (RDTM) methods. For the purpose of efficiency and accuracy, our results are compared with exact solutions (2)-(3) of eq. (1). Numerical considerations are revealed that RDTM is very effective and more convergent than DTM. In addition, RDTM can be applied easier than DTM and ensures very accurate solutions as in tables (3)-(6) and figures (1)-(4).

## 2. Basic Properties of Two Dimensional Reduced Differential Transform Method (RDTM) and Differential Transform Method (DTM)

### 2.1. Two dimensional DTM

Differential transform method (DTM) is a numerical method based on Taylor expansion. This method tries to find coefficients of series expansion of unknown function term by term. The concept of DTM was first proposed by Zhou [17]. By, considering the literature [17,18,19,20,21,22,23,24], we give the following definition of two dimensional DTM;

**Definition 2.1.:** Lets $u(x,t)$ is denoted two variables analytic and differentiable function, then two dimensional transform is follow

$$U(k,h) = \frac{1}{k!h!}\left[\frac{\partial^{k+h}}{\partial x^k \partial t^h} u(x,t)\right]_{\substack{x=x_0 \\ t=t_0}} \tag{4}$$

where $U(k,h)$ is the transformed function of $u(x,t)$. The transformation is called *T-function*. Hence, the differential inverse transform of $U(k,h)$ is defined as

$$u(x,t) = \sum_{k=0}^{\infty}\sum_{h=0}^{\infty} U(k,h)(x-x_0)^k (t-t_0)^h \tag{5}$$

From the eqs. (4) and (5), it can be written

$$u(x,t) = \sum_{k=0}^{\infty}\sum_{h=0}^{\infty} \frac{1}{k!h!}\left[\frac{\partial^{k+h}}{\partial x^k \partial t^h} u(x,t)\right]_{\substack{x=x_0 \\ t=t_0}} (x-x_0)^k (t-t_0)^h \tag{6}$$

In terms of applicability, we rearrange the eq. (6) as follow

$$u(x,t) = \sum_{k=0}^{n}\sum_{h=0}^{m} \frac{1}{k!h!}\left[\frac{\partial^{k+h}}{\partial x^k \partial t^h} u(x,t)\right]_{\substack{x=0 \\ t=0}} x^k t^h + R_{nm}(x,t) \tag{7}$$

where $(x_0, t_0)$ are taken as $(0,0)$ and

$$R_{nm}(x,t) = \sum_{k=n+1}^{\infty} \sum_{h=m+1}^{\infty} U(k,h) x^k t^h \qquad (8)$$

Here, $R_{nm}(x,t)$ is negligibly small terms. Some of the transform form of functions are given as Table 1 and their proofs can be found in [17,18,19,20].

**Table 1: Some two dimensional DTM operations with transformed forms.**

| Original functions | Transformed forms |
|---|---|
| $u(x,t) = v(x,t) \pm w(x,t)$ | $U(k,h) = V(k,h) \pm W(k,h)$ |
| $u(x,t) = \lambda v(x,t)$ | $U(k,h) = \lambda V(k,h)$ |
| $u(x,t) = \dfrac{\partial}{\partial x} v(x,t)$ | $U(k,h) = (k+1) V(k+1,h)$ |
| $u(x,t) = \dfrac{\partial}{\partial t} v(x,t)$ | $U(k,h) = (h+1) V(k,h+1)$ |
| $u(x,t) = \dfrac{\partial^{m+n}}{\partial x^m \partial t^n} v(x,t)$ | $U(k,h) = \dfrac{(k+m)!}{k!} \dfrac{(h+n)!}{h!} V(k+m, h+n)$ |
| $u(x,t) = v(x,t) w(x,t)$ | $U(k,h) = \sum_{r=0}^{k} \sum_{s=0}^{h} V(r, h-s) W(k-r, s)$ |
| $u(x,t) = v(x,t) w(x,t) q(x,t)$ | $U(k,h) = \sum_{r=0}^{k} \sum_{p=0}^{k-r} \sum_{s=0}^{h} \sum_{z=0}^{h-s} V(r, h-s-z) W(p, s) Q(k-r-p, z)$ |
| $u(x,t) = x^m t^n$ | $U(k,h) = \delta(k-m, h-n) = \delta(k-m) \delta(h-n)$  $\delta(k-m) = \begin{cases} 1, & k = m \\ 0, & otherwise \end{cases}$ |

### 2.2. Two dimensional RDTM

Reduced differential transform method (RDTM) which has an alternative approach of problems is presented to overcome the demerit complex calculation, discretization, linearization or perturbations of well-known numerical and analytical methods such as Adomian decomposition, Differential transform, Homotopy perturbation and Variational iteration. RDTM was first introduced by Keskin and Oturanc [26,27,28,29]. The main advantage of RDTM is providing an analytic approximation, in many cases exact solutions, in rapidly convergent sequence with elegantly computed terms [24],[26,27,28,29,30,31]. And also, unlike the DTM, RDTM is based on the Poisson series coefficients expansion. By using the literature [24],[26,27,28,29,30,31], we present the RDTM as follow.

**Definition 2.2.:** Lets $u(x,t)$ is two variables function and assumed that it can be demonstrated as a product of two functions which are single variable $u(x,t) = h(x)g(t)$. By making use of differential transform properties, $u(x,t)$ can be written as

$$u(x,t) = \sum_{i=0}^{\infty} H(i)x^i \sum_{j=0}^{\infty} G(j)t^j = \sum_{k=0}^{\infty} U_k(x)t^k \qquad (9)$$

Here $U_k(x)$ is called $t$ dimensional spectrum function of $u(x,t)$. If function $u(x,t)$ is analytic and differentiated continuously with respect to time $t$ and space $x$ in the domain of interest, than

$$U_k(x) = \frac{1}{k!}\left[\frac{\partial^k}{\partial t^k} u(x,t)\right]_{t=t_0} \qquad (10)$$

where $U_k(x)$ is transformed function of $u(x,t)$. The differential inverse transform of $U_k(x)$ is defined as

$$u(x,t) = \sum_{k=0}^{\infty} U_k(x)(t-t_0)^k \qquad (11)$$

Combining (9)-(11), we can write

$$u(x,t) = \sum_{k=0}^{\infty} \frac{1}{k!}\left[\frac{\partial^k}{\partial t^k} u(x,t)\right]_{t=t_0} (t-t_0)^k \qquad (12)$$

In real applications, we use the finite series form of (12), therefore we rewrite the solution as

$$\tilde{u}_n(x,t) = \sum_{k=0}^{n} U_k(x)t^k \qquad (13)$$

where $n$ is order of approximation. Hence, the RDTM solution is given by

$$u(x,t) = \lim_{n \to \infty} \tilde{u}_n(x,t) \tag{14}$$

here $n$ is taken as sufficiently big to get convergent solution. In Table 2, transformed form of mathematical operation of some functions are given and their proofs are shown in ref. [26,27].

### 3. Solution procedures of Ramani Equation by DTM and RDTM

### 3.1. DTM procedure

Let's consider the coupled Ramani equation (1) with two different initial conditions as [4],[10],[12,13,14,15,16]

$$u(x,0) = a_0 + 2\alpha \tanh(\alpha x)$$
$$v(x,0) = -\frac{4\beta\alpha^4}{9} - \frac{16\alpha^6}{27} + \frac{5\beta^2\alpha^2}{9} - \frac{5\beta^3}{54}$$
$$+ \left( \frac{20\beta\alpha^4}{9} + \frac{16\alpha^6}{9} - \frac{5\beta^2\alpha^2}{9} \right) \tanh^2(\alpha x) \tag{15}$$

$$u(x,0) = a_0 - 2\alpha \tan(\alpha x)$$
$$v(x,0) = -\frac{4\beta\alpha^4}{9} + \frac{16\alpha^6}{27} - \frac{5\beta^2\alpha^2}{9} - \frac{5\beta^3}{54}$$
$$+ \left( -\frac{20\beta\alpha^4}{9} + \frac{16\alpha^6}{9} - \frac{5\beta^2\alpha^2}{9} \right) \tan^2(\alpha x) \tag{16}$$

$U(k,h), V(k,h)$, which are called *T-function*, denote the transformation of the functions $u(x,t), v(x,t)$ in eq. (1) respectively. Then from Table 1 and (4) to (7), we obtain the transformed form of eq. (1) as below

$$\begin{cases} 5(h+1)(h+2)U(k,h+2) = \dfrac{(k+6)!}{k!}U(k+6,h) + 18(k+1)V(k+1,h) \\ +15\sum_{r=0}^{k}\sum_{s=0}^{h}(r+1)(k-r+4)(k-r+1)(k-r+2)(k-r+3) \\ \quad U(r+1,h-s)U(k-r+4,s) \end{cases} \tag{17}$$

$$\begin{cases} +15\sum_{r=0}^{k}\sum_{s=0}^{h}(r+1)(r+2)(k-r+1)(k-r+2)(k-r+3)U(r+2,h-s) \\ \phantom{+15\sum_{r=0}^{k}\sum_{s=0}^{h}}U(k-r+3,s) \\ -5\dfrac{(k+3)!}{k!}(h+1)U(k+3,h+1) \\ -15\sum_{r=0}^{k}\sum_{s=0}^{h}(h-s+1)(k-r+1)(k-r+2)U(k-r+2,s) \\ \phantom{-15\sum_{r=0}^{k}\sum_{s=0}^{h}}U(r,h-s+1) \\ -15\sum_{r=0}^{k}\sum_{s=0}^{h}(r+1)(k-r+1)(h-s+1)U(r+1,h-s)U(k-r+1,h-s+1) \\ +45\sum_{r=0}^{k}\sum_{l=0}^{k-r}\sum_{s=0}^{h}\sum_{p=0}^{h-s}(r+1)(l+1)(k-r-l+1)(k-r-l+2) \\ \phantom{+45\sum_{r=0}^{k}\sum_{l=0}^{k-r}\sum_{s=0}^{h}\sum_{p=0}^{h-s}}U(r+1,h-s-p)U(l+1,s)U(k-r-l+2,p) \end{cases}$$

$$\begin{cases} (h+1)V(k,h+1) = \dfrac{(k+3)!}{k!}V(k+3,h) \\ +3\sum_{r=0}^{k}\sum_{s=0}^{h}(r+1)(k-r+1)V(r+1,h-s) \\ \phantom{+3\sum_{r=0}^{k}\sum_{s=0}^{h}}U(k-r+1,s) \\ +3\sum_{r=0}^{k}\sum_{s=0}^{h}(k-r+1)(k-r+2)V(r,h-s) \\ \phantom{+3\sum_{r=0}^{k}\sum_{s=0}^{h}}U(k-r+2,s) \end{cases} \quad (18)$$

and for initial conditions (15-16), we obtain as

$$U(k,0) = a_0\delta(k,0) + 2\alpha\left(\dfrac{(2\alpha)^k - k!\delta(k,0)}{(2\alpha)^k + k!\delta(k,0)}\right)$$

$$V(k,0) = \left(-\dfrac{4\beta\alpha^4}{9} - \dfrac{16\alpha^6}{27} + \dfrac{5\beta^2\alpha^2}{9} - \dfrac{5\beta^3}{54}\right)\delta(k,0) \quad (19)$$

$$+ \left(\dfrac{20\beta\alpha^4}{9} + \dfrac{16\alpha^6}{9} - \dfrac{5\beta^2\alpha^2}{9}\right)\left(\dfrac{(2\alpha)^k - k!\delta(k,0)}{(2\alpha)^k + k!\delta(k,0)}\right)^2$$

and

$$U(k,0) = a_0\delta(k,0) - 2\alpha\tan\left(\dfrac{k\pi}{2}\right)$$

$$V(k,0) = \left(-\frac{4\beta\alpha^4}{9} + \frac{16\alpha^6}{27} - \frac{5\beta^2\alpha^2}{9} - \frac{5\beta^3}{54}\right)\delta(k,0) \quad (20)$$

$$+ \left(-\frac{20\beta\alpha^4}{9} + \frac{16\alpha^6}{9} - \frac{5\beta^2\alpha^2}{9}\right)\left(\tan\left(\frac{k\pi}{2}\right)\right)^2$$

We put firstly (19) into (17-18) and using the DTM procedures (4) to (7), we get the two terms approximate traveling DTM solution of coupled Ramani equation as

$$U_{2,2}(x,t) = a_0 + 2\alpha^2 x - 2\beta\alpha^2 t + 2\beta\alpha^4 x^2 t - 2xt^2\beta^2\alpha^4$$

$$V_{2,2}(x,t) = -\frac{4\beta\alpha^4}{9} - \frac{16\alpha^6}{27} + \frac{5\beta^2\alpha^2}{9} - \frac{5\beta^3}{54} + \frac{20\beta\alpha^6 x^2}{9}$$

$$+ \frac{16\alpha^8 x^2}{9} - \frac{5\beta^2\alpha^4 x^2}{9} - \frac{32\beta xt\alpha^8}{9} - \frac{40xt\beta^2\alpha^6}{9} \quad (21)$$

$$+ \frac{10xt\alpha^4\beta^3}{9} + \frac{16\beta^2\alpha^8 t^2}{9} + \frac{20\beta^3\alpha^6 t^2}{9} - \frac{5t^2\beta^4\alpha^4}{9}$$

$$- \frac{64\alpha^{10} x^2 t^2 \beta^2}{9} - \frac{80\beta^3\alpha^8 x^2 t^2}{9} + \frac{20\beta^4\alpha^6 x^2 t^2}{9}$$

and secondly put (20) into (17-18), we obtain the other traveling DTM solution of eq. (1) as following

$$U_{2,2}(x,t) = a_0 - 2\alpha^2 x + 2\beta\alpha^2 t + 2\beta\alpha^4 x^2 t - 2xt^2\beta^2\alpha^4$$

$$V_{2,2}(x,t) = -\frac{4\beta\alpha^4}{9} + \frac{16\alpha^6}{27} - \frac{5\beta^2\alpha^2}{9} - \frac{5\beta^3}{54} - \frac{20\beta\alpha^6 x^2}{9} \quad (22)$$

$$+ \frac{16\alpha^8 x^2}{9} - \frac{5\beta^2\alpha^4 x^2}{9} - \frac{32\beta xt\alpha^8}{9} + \frac{40xt\beta^2\alpha^6}{9}$$

$$+ \frac{10xt\alpha^4\beta^3}{9} + \frac{16\beta^2\alpha^8 t^2}{9} - \frac{20\beta^3\alpha^6 t^2}{9} - \frac{5t^2\beta^4\alpha^4}{9}$$

$$+ \frac{64\alpha^{10} x^2 t^2 \beta^2}{9} - \frac{80\beta^3\alpha^8 x^2 t^2}{9} - \frac{20\beta^4\alpha^6 x^2 t^2}{9}$$

Hence, it's clearly seen in Tables 3 to 6 that solutions (21) and (22) provide the good accuracy with compared exact solutions [12].

### 3.1. RDTM procedure

As the same manner, again we consider the eq. (1) with initial conditions (15-16) to obtain the RDTM solutions. $U_k(x), V_k(x)$, which are called $t$ dimensional spectrum functions, denote the transformation of the functions $u(x,t), v(x,t)$ in eq. (1) respectively. Then from Table 2 and (9) to (14), we obtain the transformed form of eq. (1) as below

$$\begin{cases} 5(k+1)(k+2)U_{k+2}(x) = \dfrac{d^6}{dx^6}U_k(x) + 18\dfrac{d}{dx}V_k(x) \\ +15\sum_{r=0}^{k}\dfrac{d^2}{dx^2}U_{k-r}(x)\dfrac{d^3}{dx^3}U_r(x) + 15\sum_{r=0}^{k}\dfrac{d}{dx}U_{k-r}(x)\dfrac{d^4}{dx^4}U_r(x) \\ +45\sum_{r=0}^{k}\sum_{s=0}^{k-r}\dfrac{d}{dx}U_r(x)\dfrac{d}{dx}U_s(x)\dfrac{d^2}{dx^2}U_{k-r-s}(x) \\ -5(k+1)\dfrac{d^3}{dx^3}U_{k+1}(x) - 15\sum_{r=0}^{k}\dfrac{d^2}{dx^2}U_r(x)(k-r+1)U_{k-r+1}(x) \\ -15\sum_{r=0}^{k}\dfrac{d}{dx}U_r(x)(k-r+1)\dfrac{d}{dx}U_{k-r+1}(x) \end{cases} \quad (23)$$

$$\begin{cases} (k+1)V_{k+1}(x) = \dfrac{d^3}{dx^3}V_k(x) \\ +3\sum_{r=0}^{k}\dfrac{d}{dx}V_r(x)\dfrac{d}{dx}U_{k-r}(x) + 3\sum_{r=0}^{k}V_r(x)\dfrac{d^2}{dx^2}U_{k-r}(x) \end{cases} \quad (24)$$

and for initial conditions (15-16), we obtain reduced transform form as respectively

$$U_0(x) = a_0 + 2\alpha \tanh(\alpha x)$$

$$V_0(x) = \left(-\dfrac{4\beta\alpha^4}{9} - \dfrac{16\alpha^6}{27} + \dfrac{5\beta^2\alpha^2}{9} - \dfrac{5\beta^3}{54}\right) + \left(\dfrac{20\beta\alpha^4}{9} + \dfrac{16\alpha^6}{9} - \dfrac{5\beta^2\alpha^2}{9}\right)\tanh(\alpha x)^2 \quad (25)$$

and

$$U_0(x) = a_0 - 2\alpha \tan(\alpha x)$$
$$V_0(x) = \left(-\frac{4\beta\alpha^4}{9} + \frac{16\alpha^6}{27} - \frac{5\beta^2\alpha^2}{9} - \frac{5\beta^3}{54}\right) \quad (26)$$
$$+ \left(-\frac{20\beta\alpha^4}{9} + \frac{16\alpha^6}{9} - \frac{5\beta^2\alpha^2}{9}\right)\tan(\alpha x)^2$$

**Table 2: Some two dimensional RDTM operations with transformed forms.**

| Original functions | Transformed forms |
|---|---|
| $u(x,t) = v(x,t) \pm w(x,t)$ | $U_k(x) = V_k(x) \pm W_k(x)$ |
| $u(x,t) = \lambda v(x,t)$ | $U_k(x) = \lambda V_k(x)$, $\lambda$ is constant |
| $u(x,t) = \dfrac{\partial}{\partial x}v(x,t)$ | $U_k(x) = \dfrac{\partial}{\partial}V_k(x)$ |
| $u(x,t) = \dfrac{\partial^r}{\partial t^r}v(x,t)$ | $U_k(x) = \dfrac{(k+r)!}{k!}V_{k+r}(x)$ |
| $u(x,t) = v(x,t)w(x,t)$ | $U_k(x) = \sum_{r=0}^{k} V_r(x)W_{k-r}(x) = \sum_{r=0}^{k} W_r(x)V_{k-r}(x)$ |
| $u(x,t) = v(x,t)w(x,t)q(x,t)$ | $U_k(x) = \sum_{r=0}^{k}\sum_{p=0}^{k-r} V_r(x)W_p(x)Q_{k-r-p}(x)$ |
| $u(x,t) = x^m t^n$ | $U_k(x) = x^m \delta(k-n)$, $\delta(k-n) = \begin{cases} x^m, & k=n \\ 0, & \text{otherwise} \end{cases}$ |

As in the DTM solution process, by using RDTM algorithm we put firstly (25) into (23-24), we get the two terms RDTM solution of coupled Ramani equation as

$$U_2(x) = \frac{1}{\cosh(\alpha x)^7}\begin{pmatrix} a_0\cosh(\alpha x)^7 + 2\alpha\sinh(\alpha x)\cosh(\alpha x)^6 \\ -2\beta\alpha^2 t\cosh(\alpha x)^5 + 96\alpha^7 t^2\sinh(\alpha x)\cosh(\alpha x)^2 \\ -144t^2\alpha^7\sinh(\alpha x) + 24t^2\alpha^5\beta\sinh(\alpha x)\cosh(\alpha x)^2 \\ -2t^2\alpha^3\beta^2\sinh(\alpha x)\cosh(\alpha x)^4 \end{pmatrix} \quad (27)$$

$$V_2(x) = \frac{1}{54\cosh(\alpha x)^8}\begin{pmatrix} 302400t^2\alpha^{10}\beta - 75600t^2\alpha^8\beta^2 \\ +96\beta\alpha^4\cosh(\alpha x)^8 + 241920t^2\alpha^{12} \end{pmatrix} \quad (28)$$

$$\frac{1}{54\cosh(\alpha x)^8} \begin{pmatrix} +64\alpha^6 \cosh(\alpha x)^8 - 5\beta^3 \cosh(\alpha x)^8 \\ +960t\alpha^7 \beta \sinh(\alpha x)\cosh(\alpha x)^5 \\ -2880t\alpha^7 \beta \sinh(\alpha x)\cosh(\alpha x)^3 \\ -240t\alpha^5 \beta^2 \sinh(\alpha x)\cosh(\alpha x)^5 \\ +720t\alpha^5 \beta^2 \sinh(\alpha x)\cosh(\alpha x)^3 \\ -t^2\alpha^{10}\beta\left(3840\cosh(\alpha x)^6 - 120960\cosh(\alpha x)^4\right) \\ -403200t^2\alpha^{10}\beta \cosh(\alpha x)^2 + 960t^2\alpha^8\beta^2 \cosh(\alpha x)^6 \\ -30240t^2\alpha^8\beta^2 \cosh(\alpha x)^4 + 100800t^2\alpha^8\beta^2 \cosh(\alpha x)^2 \\ -120\alpha^4\beta \cosh(\alpha x)^6 + 30\alpha^2\beta^2 \cosh(\alpha x)^6 \\ -3072\alpha^{12}t^2 \cosh(\alpha x)^6 + 96768\alpha^{12}t^2 \cosh(\alpha x)^4 \\ -322560\alpha^{12}t^2 \cosh(\alpha x)^2 + 768\alpha^9 t \cosh(\alpha x)^5 \sinh(\alpha x) \\ -2304\alpha^9 t \cosh(\alpha x)^3 \sinh(\alpha x) - 96\alpha^6 \cosh(\alpha x)^6 \end{pmatrix}$$

and secondly put (26) into (23-24), we obtain the other RDTM solution of eq. (1) as following

$$U_2(x) = -\frac{1}{\cosh(\alpha x)^7} \begin{pmatrix} -a_0 \cos(\alpha x)^7 + 2\alpha \sin(\alpha x)\cos(\alpha x)^6 \\ -2\beta\alpha^2 t \cos(\alpha x)^5 - 96\alpha^7 t^2 \sin(\alpha x)\cos(\alpha x)^2 \\ +144t^2\alpha^7 \sin(\alpha x) + 24t^2\alpha^5\beta \sin(\alpha x)\cos(\alpha x)^2 \\ +2t^2\alpha^3\beta^2 \sin(\alpha x)\cos(\alpha x)^4 \end{pmatrix} \quad (29)$$

$$V_2(x) = -\frac{1}{54\cosh(\alpha x)^8} \begin{pmatrix} 302400t^2\alpha^{10}\beta + 75600t^2\alpha^8\beta^2 \\ -96\beta\alpha^4 \cos(\alpha x)^8 - 241920t^2\alpha^{12} \\ +64\alpha^6 \cos(\alpha x)^8 + 5\beta^3 \cos(\alpha x)^8 \\ -960t\alpha^7 \beta \sin(\alpha x)\cos(\alpha x)^5 \\ +2880t\alpha^7 \beta \sin(\alpha x)\cos(\alpha x)^3 \\ -240t\alpha^5 \beta^2 \sin(\alpha x)\cos(\alpha x)^5 \\ +720t\alpha^5 \beta^2 \sin(\alpha x)\cos(\alpha x)^3 \\ -t^2\alpha^{10}\beta\left(3840\cos(\alpha x)^6 - 120960\cos(\alpha x)^4\right) \\ +120\alpha^4\beta \cos(\alpha x)^6 + 30\alpha^2\beta^2 \cos(\alpha x)^6 \end{pmatrix} \quad (30)$$

$$-\frac{1}{54\cosh(\alpha x)^8}\begin{pmatrix} -403200t^2\alpha^{10}\beta\cos(\alpha x)^2 - 960t^2\alpha^8\beta^2\cos(\alpha x)^6 \\ +30240t^2\alpha^8\beta^2\cos(\alpha x)^4 - 100800t^2\alpha^8\beta^2\cos(\alpha x)^2 \\ +322560\alpha^{12}t^2\cos(\alpha x)^2 + 768\alpha^9 t\cos(\alpha x)^5\sin(\alpha x) \\ +3072\alpha^{12}t^2\cos(\alpha x)^6 - 96768\alpha^{12}t^2\cos(\alpha x)^4 \\ -2304\alpha^9 t\cos(\alpha x)^3\sin(\alpha x) - 96\alpha^6\cos(\alpha x)^6 \end{pmatrix}$$

Thus, it's obviously noted on the Tables 3 to 6 and Figures 1 to 4 that solutions (29) and (30) provide the good accuracy with compared exact [12] and DTM solutions.

**Table 3:** The numerical results of three terms DTM and RDTM solutions of eq. (1) with compared exact solution (2) at $t = 20$ and $a_0 = 1, \alpha = \beta = 0.01$

| $x$ | Exact [12] $u(x,t)$ | RDTM $U_3(x)$ | DTM $U_{3,3}(x,t)$ | Exact [12] $v(x,t)$ | RDTM $V_3(x)$ | DTM $V_{3,3}(x,t)$ |
|---|---|---|---|---|---|---|
| -50 | 0.990726228 | 0.9907262253 | 0.9908033667 | $-8.822839648\times10^{-8}$ | $-8.822106112\times10^{-8}$ | $-8.842207168\times10^{-8}$ |
| -40 | 0.9923668212 | 0.9923668187 | 0.9923930999 | $-8.785868735\times10^{-8}$ | $-8.785217835\times10^{-8}$ | $-8.794184127\times10^{-8}$ |
| -30 | 0.9941371636 | 0.9941371613 | 0.9941436287 | $-8.754022597\times10^{-8}$ | $-8.753492301\times10^{-8}$ | $-8.756755781\times10^{-8}$ |
| -20 | 0.9960140671 | 0.9960140653 | 0.9960149542 | $-8.729383894\times10^{-8}$ | $-8.729008043\times10^{-8}$ | $-8.729939194\times10^{-8}$ |
| -10 | 0.9979670454 | 0.9979670444 | 0.9979670778 | $-8.713716109\times10^{-8}$ | $-8.713520228\times10^{-8}$ | $-8.712751428\times10^{-8}$ |
| 10 | 1.001953749 | 1.001953750 | 1.001953722 | $-8.713295226\times10^{-8}$ | $-8.713487060\times10^{-8}$ | $-8.713330592\times10^{-8}$ |
| 20 | 1.003909050 | 1.003909050 | 1.003908246 | $-8.72857484\times10^{-8}$ | $-8.728947099\times10^{-8}$ | $-8.729131646\times10^{-8}$ |
| 30 | 1.005789625 | 1.005789629 | 1.005783571 | $-8.752885527\times10^{-8}$ | $-8.753412915\times10^{-8}$ | $-8.755629763\times10^{-8}$ |
| 40 | 1.007564728 | 1.007564730 | 1.0075397 | $-8.784482114\times10^{-8}$ | $-8.785130927\times10^{-8}$ | $-8.792842001\times10^{-8}$ |
| 50 | 1.009210856 | 1.009210859 | 1.009136633 | $-8.821289529\times10^{-8}$ | $-8.822021826\times10^{-8}$ | $-8.840785424\times10^{-8}$ |

## 4. Results and Conclusions

In this paper, we consider the very famous physical problems coupled Ramani equation (1) to find two approximate traveling wave solutions by using DTM and RDTM. Moreover, we perfectly obtain approximate solutions of (1) compatible with exact solutions in [12]. In order to test efficiency, convergence and accuracy of DTM and RDTM, we perform the numerical values $a_0 = 1, \alpha = 0.01, \beta = 0.01$ in the three term approximate solutions of eq. (1) which are shown in Tables 3 to 6. Also, for $a_0 = \frac{1}{2}, \alpha = 0.03, \beta = 0.03$ and $a_0 = \frac{1}{2}, \alpha = 0.02, \beta = 0.02$, error rates and comparisons of exact [12] and three terms RDTM, DTM solutions are presented in Figures 1 to 4. Results show that DTM and RDTM are efficient and powerful technique where RDTM is more easier, fast and better than DTM.

**Table 4:** Absolute errors between exact (2) and three terms DTM and RDTM solutions of eg. (1) at $t = 20$ and $a_0 = 1, \alpha = \beta = 0.01$

| $x$ | RDTM $\lvert u(x,t) - U_3(x) \rvert$ | DTM $\lvert u(x,t) - U_{3,3}(x,t) \rvert$ | RDTM $\lvert v(x,t) - V_3(x) \rvert$ | DTM $\lvert v(x,t) - V_{3,3}(x,t) \rvert$ |
|---|---|---|---|---|
| -50 | $2.7 \times 10^{-9}$ | $7.71387 \times 10^{-5}$ | $7.33534 \times 10^{-12}$ | $1.9367522 \times 10^{-10}$ |
| -40 | $2.5 \times 10^{-9}$ | $2.62787 \times 10^{-5}$ | $6.509 \times 10^{-12}$ | $8.315392 \times 10^{-11}$ |
| -30 | $2.3 \times 10^{-9}$ | $6.4651 \times 10^{-6}$ | $5.30296 \times 10^{-12}$ | $2.733184 \times 10^{-11}$ |
| -20 | $1.8 \times 10^{-9}$ | $8.871 \times 10^{-7}$ | $3.75851 \times 10^{-12}$ | $5.553 \times 10^{-12}$ |
| -10 | $1.0 \times 10^{-9}$ | $3.24 \times 10^{-8}$ | $1.95881 \times 10^{-12}$ | $3.5319 \times 10^{-13}$ |
| 10 | $1.0 \times 10^{-9}$ | $2.7 \times 10^{-8}$ | $1.91834 \times 10^{-12}$ | $3.5366 \times 10^{-13}$ |
| 20 | 0.0 | $8.04 \times 10^{-7}$ | $3.72259 \times 10^{-12}$ | $5.56806 \times 10^{-12}$ |
| 30 | $4.0 \times 10^{-9}$ | $6.054 \times 10^{-6}$ | $5.27388 \times 10^{-12}$ | $2.744236 \times 10^{-11}$ |
| 40 | $2.0 \times 10^{-9}$ | $2.5028 \times 10^{-5}$ | $6.48813 \times 10^{-12}$ | $8.359887 \times 10^{-11}$ |
| 50 | $3.0 \times 10^{-9}$ | $7.4223 \times 10^{-5}$ | $7.32297 \times 10^{-12}$ | $1.9495895 \times 10^{-10}$ |

**Table 5:** The numerical results of three terms DTM and RDTM solutions of eq. (1) with compared exact solution (3) at $t = 20$ and $a_0 = 1, \alpha = \beta = 0.01$

| $x$ | Exact [12] $u(x,t)$ | RDTM $U_3(x)$ | DTM $U_{3,3}(x,t)$ | Exact [12] $v(x,t)$ | RDTM $V_3(x)$ | DTM $V_{3,3}(x,t)$ |
|---|---|---|---|---|---|---|
| -50 | 1.010978045 | 1.010978054 | 1.010883361 | $-9.993227215 \times 10^{-8}$ | $-9.991393179 \times 10^{-8}$ | $-9.965144899 \times 10^{-8}$ |
| -40 | 1.008503055 | 1.008503059 | 1.008473089 | $-9.923603996 \times 10^{-8}$ | $-9.922331726 \times 10^{-8}$ | $-9.91274111 \times 10^{-8}$ |
| -30 | 1.006230580 | 1.006230583 | 1.006223615 | $-9.875256259 \times 10^{-8}$ | $-9.874398368 \times 10^{-8}$ | $-9.871963441 \times 10^{-8}$ |
| -20 | 1.004095861 | 1.004095864 | 1.004094943 | $-9.843424659 \times 10^{-8}$ | $-9.842893010 \times 10^{-8}$ | $-9.842793407 \times 10^{-8}$ |
| -10 | 1.002047104 | 1.002047105 | 1.002047072 | $-9.825251278 \times 10^{-8}$ | $-9.824995748 \times 10^{-8}$ | $-9.825212525 \times 10^{-8}$ |
| 10 | 0.9980337012 | 0.9980337002 | 0.9980337285 | $-9.824782982 \times 10^{-8}$ | $-9.825033776 \times 10^{-8}$ | $-9.824744281 \times 10^{-8}$ |
| 20 | 0.9959874261 | 0.9959874240 | 0.9959882569 | $-9.842449482 \times 10^{-8}$ | $-9.842975832 \times 10^{-8}$ | $-9.841819953 \times 10^{-8}$ |
| 30 | 0.9938570755 | 0.9938570718 | 0.9938635848 | $-9.873690092 \times 10^{-8}$ | $-9.874541644 \times 10^{-8}$ | $-9.870410843 \times 10^{-8}$ |
| 40 | 0.9915912460 | 0.99159124 | 0.9916197122 | $-9.921301115 \times 10^{-8}$ | $-9.922565352 \times 10^{-8}$ | $-9.910498468 \times 10^{-8}$ |
| 50 | 0.9891258314 | 0.9891258221 | 0.9892166388 | $-9.98994945 \times 10^{-8}$ | $-9.991772791 \times 10^{-8}$ | $-9.962064343 \times 10^{-8}$ |

**Table 6:** Absolute errors between exact (3) and three terms DTM and RDTM solutions of eg. (1) at $t = 20$ and $a_0 = 1, \alpha = \beta = 0.01$

| $x$ | RDTM $|u(x,t) - U_3(x)|$ | DTM $|u(x,t) - U_{3,3}(x,t)|$ | RDTM $|v(x,t) - V_3(x)|$ | DTM $|v(x,t) - V_{3,3}(x,t)|$ |
|---|---|---|---|---|
| -50 | $9.0 \times 10^{-9}$ | $9.4684 \times 10^{-5}$ | $1.834036 \times 10^{-11}$ | $2.8082316 \times 10^{-10}$ |
| -40 | $4.0 \times 10^{-9}$ | $2.9966 \times 10^{-5}$ | $1.272270 \times 10^{-11}$ | $1.0862886 \times 10^{-10}$ |
| -30 | $3.0 \times 10^{-9}$ | $6.965 \times 10^{-6}$ | $8.57891 \times 10^{-12}$ | $3.292818 \times 10^{-11}$ |
| -20 | $3.0 \times 10^{-9}$ | $9.18 \times 10^{-7}$ | $5.31649 \times 10^{-12}$ | $6.31252 \times 10^{-12}$ |
| -10 | $1.0 \times 10^{-9}$ | $3.2 \times 10^{-8}$ | $2.55530 \times 10^{-12}$ | $3.8753 \times 10^{-13}$ |
| 10 | $1.0 \times 10^{-9}$ | $2.73 \times 10^{-8}$ | $2.50794 \times 10^{-12}$ | $3.8701 \times 10^{-13}$ |
| 20 | $2.1 \times 10^{-9}$ | $8.308 \times 10^{-7}$ | $5.2635 \times 10^{-12}$ | $6.29529 \times 10^{-12}$ |
| 30 | $3.7 \times 10^{-9}$ | $6.5093 \times 10^{-6}$ | $8.51552 \times 10^{-12}$ | $3.279249 \times 10^{-11}$ |
| 40 | $6.0 \times 10^{-9}$ | $2.84662 \times 10^{-5}$ | $1.264237 \times 10^{-11}$ | $1.0802647 \times 10^{-10}$ |
| 50 | $9.3 \times 10^{-9}$ | $9.08074 \times 10^{-5}$ | $1.823341 \times 10^{-11}$ | $2.7885107 \times 10^{-10}$ |

R E F E R E N C E S

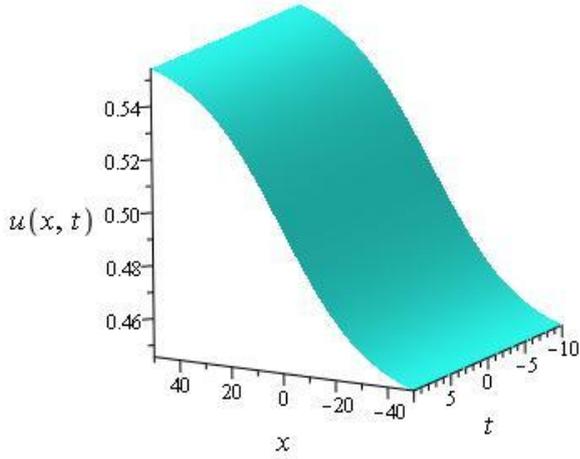

a) Exact solution of $u(x, t)$

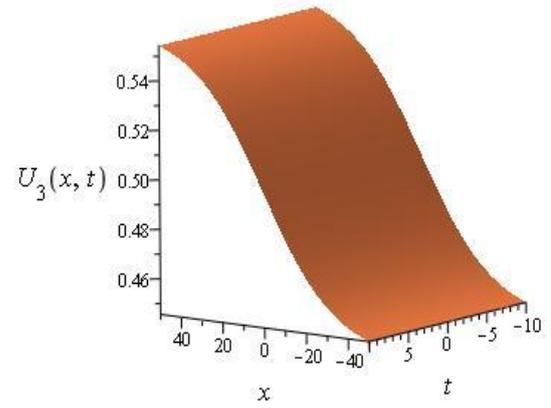

b) RDTM solution of $u(x, t)$

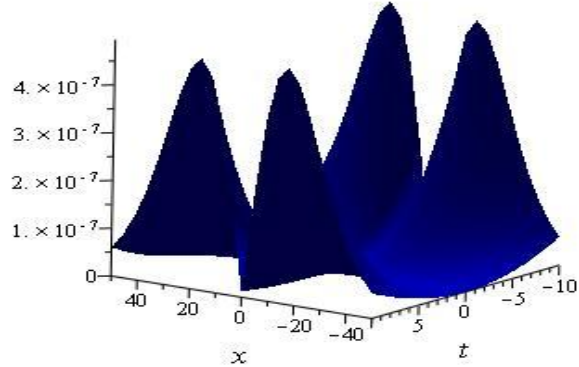

c) Errors of $|u(x, t) - U_3(x, t)|$

**Figure 1:** Comparison between exact solution (2) and three terms RDTM solution of $u(x,t)$ with errors for $a_0 = \dfrac{1}{2}, \alpha = \beta = 0.03$.

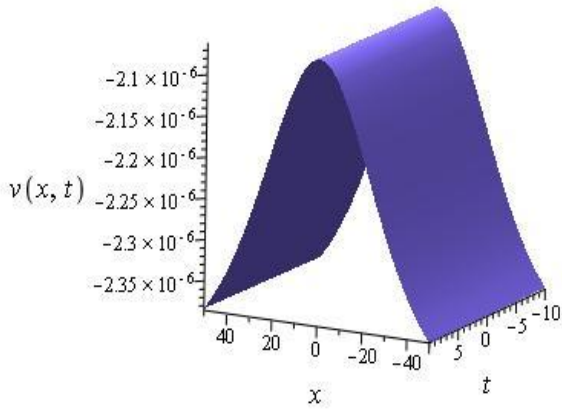

a) Exact solution of $v(x, t)$

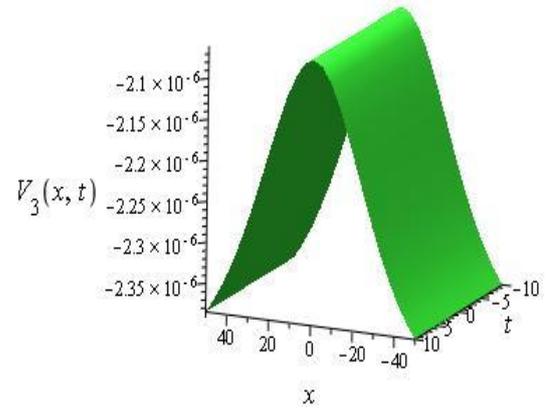

b) RDTM solution of $v(x, t)$

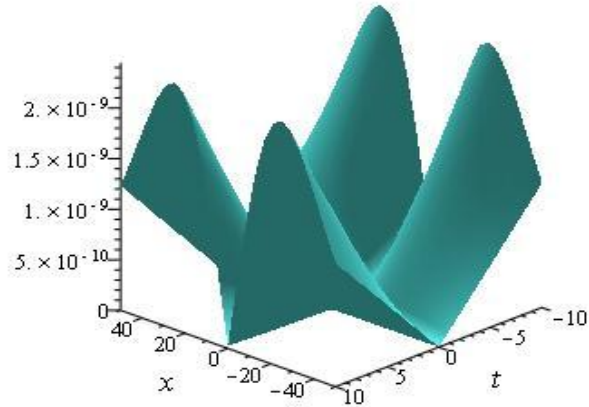

c) Errors of $|v(x, t) - V_3(x, t)|$

**Figure 2:** Comparison between exact solution (2) and three terms RDTM solution of $v(x,t)$ with errors for $a_0 = \dfrac{1}{2}, \alpha = \beta = 0.03$.

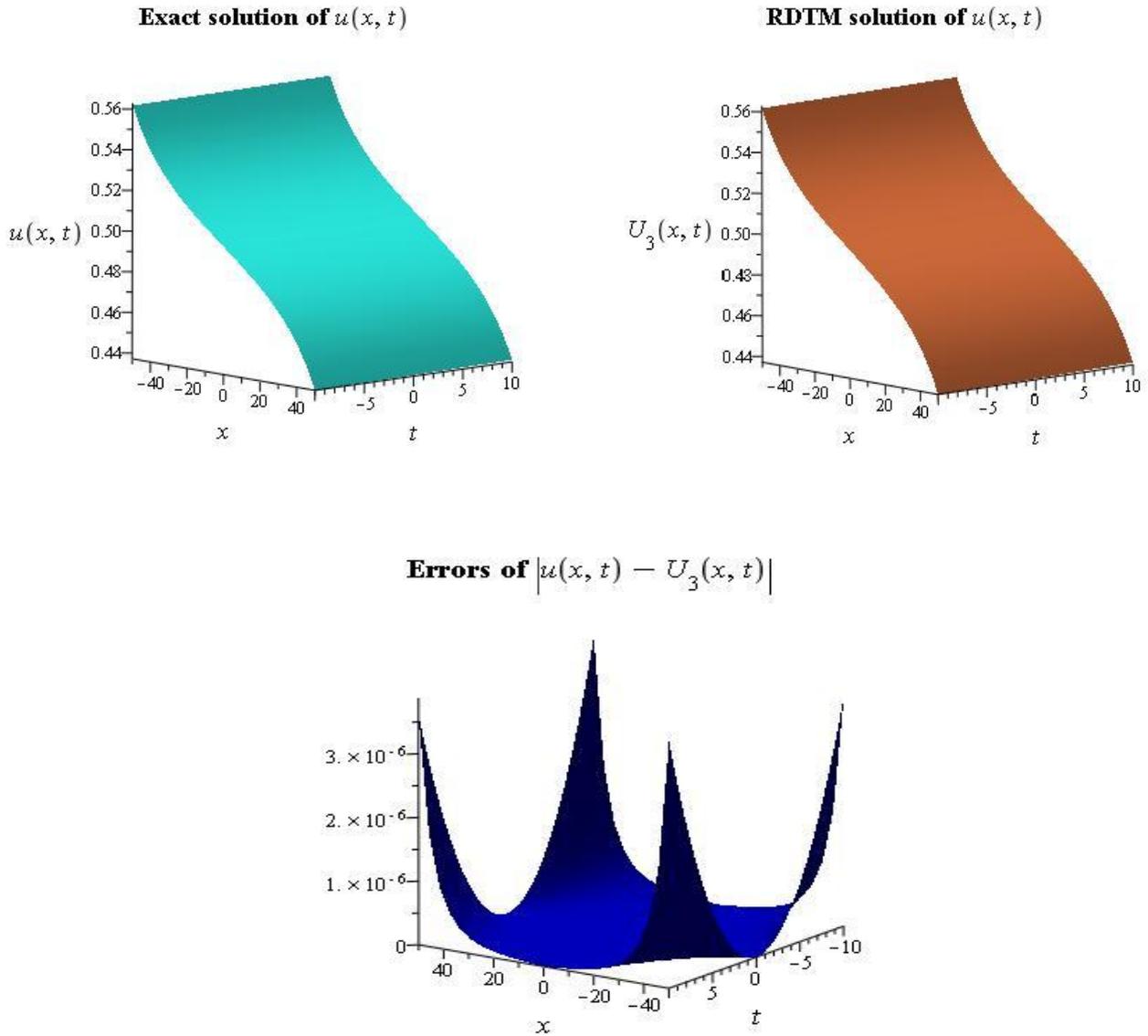

**Figure 3:** Comparison between exact solution (3) and three terms RDTM solution of $u(x,t)$ with errors for $a_0 = \dfrac{1}{2}, \alpha = \beta = 0.02$.

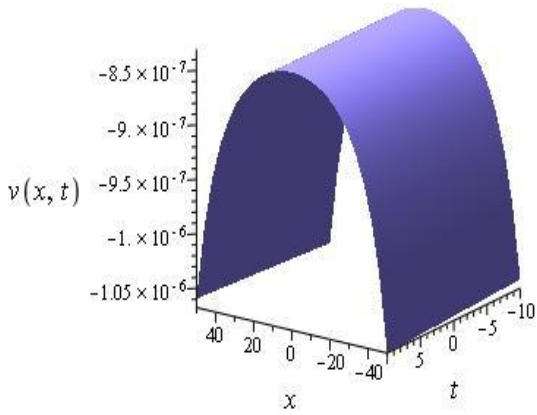
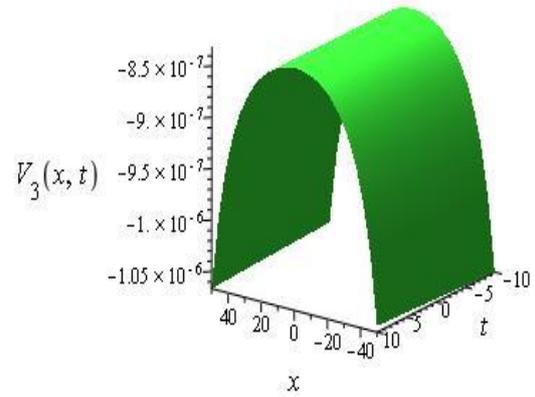
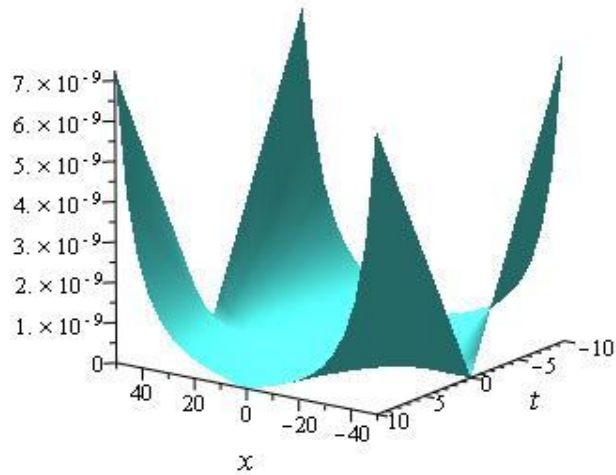

**Figure 4:** Comparison between exact solution (3) and three terms RDTM solution of $v(x,t)$ with errors for $a_0 = \dfrac{1}{2}, \alpha = \beta = 0.02$.